\documentclass[11pt]{amsart}
\usepackage{amsfonts}
\usepackage{amssymb}
\usepackage{amsmath}
\usepackage{amsthm}
\usepackage{epic}
\usepackage{verbatim}
\usepackage{a4wide}
\usepackage{latexsym}
\usepackage{verbatim}
\usepackage{url}
\usepackage{graphicx,color}  
\usepackage{mathrsfs}
\def\t{\mbox{\textbf{\textsf{t}}}}

\newcommand{\NN}{\mathbb{N}}

\newcommand{\ZZ}{\mathbb{Z}}
\newcommand{\RR}{\mathbb{R}}
\newcommand{\TT}{\mathbb{T}}

\newcommand{\Pref}{\mbox{Pref}}
\newcommand{\Orbit}{\mathcal{O}}
\newcommand{\Closure}{\overline{\Orbit}}
\def\ee{\'{e}\'{e}}

\def\t1{\boldsymbol{t}^{(1)}}
\def\bt{\boldsymbol{t}}

\def\m1{\boldsymbol{m}^{(1)}}

\def\bs{\boldsymbol{s}}
\def\cA{\mathcal{A}}
\def\cS{\mathcal{S}}

\def\bbf{\boldsymbol{f}}

\def\br{\boldsymbol{r}}
\def\bu{\boldsymbol{u}}
\def\bv{\boldsymbol{v}}
\def\bx{\boldsymbol{x}}
\def\by{\boldsymbol{y}}

\def\bc{\boldsymbol{c}}
\def\rev{\widetilde}

\def\cAw{\mathcal{A}^\omega}

\def\cL{\mathcal{L}}

\def\empt{\varepsilon}

\def \proof{\medbreak\noindent{\it Proof.\ \ }}

\def \endpf{{\ \ $\Box$ \smallbreak}}

\theoremstyle{plain}
\newtheorem{theorem}{Theorem} 
\newtheorem{lemma}[theorem]{Lemma}
\newtheorem{corollary}[theorem]{Corollary}
\newtheorem{proposition}[theorem]{Proposition}

\theoremstyle{definition}
\newtheorem{definition}[theorem]{Definition}
\newtheorem{example}[theorem]{Example}

\theoremstyle{remark}
\newtheorem{remark}[theorem]{Remark}
\newtheorem*{note}{Note}

\begin{document}

\title{Extremal properties of (epi)Sturmian sequences and 
distribution modulo $1$}

\date{July 14, 2009}

\author{Jean-Paul Allouche}

\address{(J.-P.~Allouche) CNRS, LRI, UMR 8623, Universit\'e Paris-Sud, B\^atiment 490, 
F-91405 Orsay Cedex, FRANCE}

\email{allouche@lri.fr}

\author{Amy Glen}

\curraddr{(A.~Glen) Department of Mathematics and Statistics, School of Chemical and Mathematical 
Sciences, Murdoch University, Perth, WA 6150, AUSTRALIA}

\address{The Mathematics Institute, Reykjav\'ik University, Kringlan 1, 
IS-103 Reykjav\'ik, ICELAND}

\email{amy.glen@gmail.com}

\subjclass[2000]{11J71; 68R15; 11B85; 37B10}

\keywords{combinatorics on words; lexicographic order; Sturmian sequence;
episturmian sequence; Arnoux-Rauzy sequence; distribution modulo $1$; Mahler $Z$-number}

\begin{abstract}
Starting from a study of Y.\ Bugeaud and A.\ Dubickas (2005) on a question 
in distribution of real numbers modulo $1$ {\em via\,} combinatorics on words, 
we survey some combinatorial properties of (epi)Sturmian sequences and
distribution modulo $1$ in connection to their work. In particular we focus
on extremal properties of (epi)Sturmian sequences, some of which have been
rediscovered several times.
\end{abstract}

\maketitle


\section{Introduction}

A little while ago, JPA came across a paper of Y. Bugeaud and 
A. Dubickas \cite{yBaD05frac} where the authors describe all 
irrational numbers $\xi > 0$ such that the fractional parts 
$\{\xi b^n\}$, $n \geq 0$, all belong to an interval of length 
$1/b$, where $b \geq 2$ is a given integer. They also prove
that $1/b$ is the minimal length having this property. An 
interesting and unexpected result in their paper is the following:
{\em the irrational numbers $\xi > 0$ such that the fractional parts
$\{\xi b^n\}$, $n \geq 0$, all belong to a closed interval of length 
$1/b$ are exactly the positive real numbers whose base $b$ expansions 
are characteristic Sturmian sequences on $\{k, k+1\}$, where 
$k \in \{0, 1, \ldots, b-2\}$}.
(Recall that characteristic Sturmian sequences are codings of 
trajectories on a square billiard that start from a corner with 
an irrational slope; alternatively a characteristic Sturmian
sequence can be obtained by coding the sequence of cuts in an 
integer lattice over the positive quadrant of $\RR^2$ made by 
a line of irrational slope through the origin.) We will see that 
the combinatorial results underlying \cite{yBaD05frac} were stated 
several times, in particular by P.~Veerman who proved
Bugeaud-Dubickas' number-theoretical statement in the case $b=2$ 
as soon as 1986--1987 (see \cite{pV86symb, pV87symb}).

\section{The combinatorial background of a result of Bugeaud and Dubickas}

The main result of Bugeaud and Dubickas \cite[Theorem~2.1]{yBaD05frac}
will be recalled in Section~\ref{morsesection}. Looking at the proof, 
we see that its core is a result in combinatorics on
words that is encompassed by Theorems~\ref{P:JPA} and \ref{T:intro} below. 

\subsection{Sturmian sequences show up}

In this section sequences take their values in $\{0,1\}$. We let $T$ denote
the shift map defined as follows: if $\bs := (s_n)_{n \geq 0}$, then 
$T(\bs) = T((s_n)_{n \geq 0}) := (s_{n+1})_{n \geq 0}$, and we let $\leq$ 
denote the lexicographical order on $\{0,1\}^{\mathbb N}$ induced by $0 < 1$.

\begin{theorem} \label{P:JPA} 
An aperiodic sequence $\bs := (s_n)_{n \geq 0}$ on $\{0,1\}$ is Sturmian if 
and only if there exists a sequence $\bu := (u_n)_{n \geq 0}$ on $\{0,1\}$ 
such that $0\bu \leq T^k(\bs) \leq 1\bu$ for all $k \geq 0$. Moreover, $\bu$ is 
the unique characteristic Sturmian sequence with the same slope as $\bs$, and we 
have $0 \bu = \inf\{T^k(\bs), \ k \geq 0\}$ and $1 \bu = \sup\{T^k(\bs), \ k \geq 0\}$. 
\end{theorem}

\begin{theorem} \label{T:intro}
An aperiodic sequence $\bu$ on $\{0,1\}$ is a characteristic 
Sturmian sequence if and only if, for all $k \geq 0$,
$$
0\bu < T^k(\bu) < 1\bu.
$$
Furthermore, we have $0 \bu = \inf\{T^k(\bu), \ k \geq 0\}$ and 
$1 \bu = \sup\{T^k(\bu), \ k \geq 0\}$.
\end{theorem}

\noindent
[Theorem~\ref{T:intro} is an easy consequence of Theorem~\ref{P:JPA}.
For a proof of Theorem~\ref{P:JPA}, see Section~\ref{SS:Pirillo}.]

\bigskip

Actually Theorem~\ref{T:intro}  was known prior to \cite{yBaD05frac}. 
It was indicated to JPA by G. Pirillo (who published it in \cite{gP03ineq}): 
JPA suggested that this could well be already in a paper by S.~Gan 
\cite{sG01stur} under a slightly disguised form (which is indeed the case).
About 8 years earlier J. Berstel and P. S\'e\'ebold \cite{jBpS93acha} and 
also J.-P. Borel and F. Laubie \cite{jBfL93quel} proved one direction of 
Theorem~\ref{T:intro}, namely that characteristic Sturmian sequences satisfy the 
inequalities $0\bu < T^k(\bu) < 1\bu$ for all $k\geq 0$. In fact, it seems that 
both theorems were proved for the first time (including the number-theoretical 
aspect for the case of base $2$) by P. Veerman \cite{pV86symb, pV87symb}. For 
more on the history of that result (including other papers like \cite{sBpS94orde}), 
see Section~\ref{SS:Pirillo} (in particular Section~\ref{veer}).

\subsection{Generalizations}

Two directions for generalizations are possible. One is purely
combinatorial and looks at generalizations of Sturmian sequences; in particular  
{\it episturmian sequences}, which share many properties with Sturmian sequences and 
have similar extremal properties. In this direction, characterizations of finite and 
infinite (epi)Sturmian sequences via lexicographic orderings have recently been studied (see 
\cite{aG06acha, aG07orde, aGjJgP06char, oJlZ04char, jJgP02onac, gP03ineq, gP05ineq, gP05mors}).
The other type of generalization is number-theoretic and looks at distribution modulo $1$ from 
a combinatorial point of view. Recent papers of Dubickas go in this direction; we cite two
of them showing an unexpected occurrence of the Thue-Morse sequence \cite{aD06onth, aD07onas} 
(see Section~\ref{morsesection}).

\section{More on Sturmian and episturmian sequences}

We give in this section some background on Sturmian and episturmian sequences.

\subsection{Terminology \& notation}

In what follows, we shall use the following terminology and notation from combinatorics
on words (see, e.g., \cite{mL02alge}).

Let $\cA$ denote a finite non-empty {\it alphabet}. 
If $w = x_{1}x_{2}\cdots x_{m}$ is a {\em finite word}
over $\cA$, where each $x_i \in \cA$, then the {\em length} of $w$ is $|w| := m$, and we let 
$|w|_a$ denote the number of occurrences of a letter $a$ in $w$. The word of length 0 is called 
the {\em empty word}, denoted by $\empt$. The {\em reversal} $\rev{w}$ of $w$ is given by 
$\rev{w} = x_{m}x_{m-1}\cdots x_{1}$, and if $w = \rev{w}$, then $w$ is called a 
{\em palindrome}. 

An {\em infinite word} (or simply {\em sequence}) $\bx$ over $\cA$ is a sequence indexed by 
$\NN$ with values in $\cA$, i.e., $\bx = x_0x_1x_2\cdots$, where each $x_i \in \cA$. A finite 
word $w$ is a {\em factor} of $\bx$ if $w = \varepsilon$ or 
$w = x_i\cdots x_{j}$ for some $i$, $j$ with $i \leq j$. Furthermore, if $w$ is not empty, 
$w$ is said to be a {\em prefix} of $\bx$ if $i=0$, and we say that $w$ is 
{\em right} (resp.~{\em left}) {\em special} if $wa$, $wb$ (resp.~$aw$, $bw$) are factors 
of $\bx$ for some letters $a$, $b \in \cA$, $a \ne b$. The set of all factors of $\bx$ is 
denoted by $F(\bx)$, and $F_n(\bx)$ denotes the set of factors of length $n$ of $\bx$, i.e., 
$F_n(\bx) := F(\bx) \cap \cA^n$. Moreover, the {\em alphabet} of $\bx$ is Alph$(\bx)
:= F(\bx) \cap \cA$. A factor of an infinite word $\bx$ is {\em recurrent} in $\bx$ if it 
occurs infinitely many times in $\bx$. The sequence $\bx$ itself is said to be 
{\em recurrent} if all of its factors are recurrent in it. Moreover $\bx$ is said 
to be {\em uniformly recurrent} (or {\em minimal}) if it is recurrent 
and if, for any factor, the gaps between its consecutive occurrences are bounded.

If $u$, $v$ are non-empty words over $\cA$, then $v^\omega$ (resp.~$uv^\omega$) denotes the 
{\it periodic} (resp.~{\it ultimately periodic}) infinite word $vvv\cdots$ (resp.~$uvvv\cdots$) 
having $|v|$ as a {\it period}. An infinite word that is not ultimately periodic is said to
be {\em aperiodic}.

For any infinite word $\bx = x_0x_1x_2x_3\cdots$, recall that the {\em shift map} $T$ is 
defined by $T(\bx) = x_1x_2x_3\cdots$. This operator naturally extends to finite words as 
a {\em circular shift} by defining $T(xw) = wx$ for any letter $x$ and finite word $w$.

The set of all finite (resp.~infinite) words over $\cA$ is denoted by $\cA^*$ 
(resp.~$\cAw$), and we define $\cA^+ := \cA^*\setminus\{\empt\}$, the set of all 
{\em non-empty} words over $\cA$. 

\subsection{Sturmian sequences} \label{SS:Sturmian}

Sturmian sequences were introduced in \cite{gHmM40symb}. They are in some sense the 
``least complicated'' aperiodic sequences on a binary alphabet, as is evident from
Lemma~\ref{easy} and Theorem~\ref{stu} below. The following lemma can essentially be 
found in \cite{gHmM40symb}.

\begin{lemma} {\rm \cite{gHmM40symb}} \label{easy}
Let $\bs$ be a sequence taking exactly $a \geq 2$ distinct values.
Let $p(k)$ be the number of distinct factors of length $k$ of $\bs$
(the function $k \mapsto p(k)$ is called the {\em block-complexity} of
the sequence $\bs$).
Then the following properties are equivalent.
\begin{itemize}
\item[(i)]   There exists $k_0 \geq 1$ such that $p(k_0+1) = p(k_0)$.
\item[(ii)]  The sequence $(p(k))_{k \geq 1}$ is ultimately constant 
(i.e., constant from some index on).
\item[(iii)] There exists $M$ such that $p(k) \leq M$ for all $k \geq 1$.
\item[(iv)]  There exists $k_1 \geq 1$ such that $p(k_1) \leq k_1 + a - 2$.
\item[(v)]   Let $g(k) = p(k) - k$. There exists $k_2 \geq 1$ such that 
             $g(k_2+1) < g(k_2)$. 
\item[(vi)]  The sequence $\bs$ is ultimately periodic.
\end{itemize}
\end{lemma}
\proof
For any sequence, we clearly have $p(k+1) \geq p(k)$ for all $k \geq 0$. 
This implies on the one hand that properties (ii) and 
(iii) are equivalent. On the other hand, this implies the 
equivalence of properties (i) and (v). Namely letting 
$g(k) := p(k) - k$, we have $g(k+1) - g(k) = p(k+1) - p(k) - 1$.

The implications (vi) $\Rightarrow$ (ii) $\Rightarrow$ (iv)
are straightforward. It thus suffices to prove that 
(iv) $\Rightarrow$ (i), and (i) $\Rightarrow$ (vi).

(iv) $\Rightarrow$ (i): if (i) is not true, then the
sequence $(p(k))_{k \geq 0}$ is (strictly) increasing. Thus, for all 
$k \geq 1$, one has $p(k+1) \geq p(k)+1$. Hence, by an easy induction,
one has $p(k) \geq p(1) + k - 1 = a + k - 1$, i.e., 
$p(k) > a + k - 2$, for all $k \geq 1$.

(i) $\Rightarrow$ (vi): the equality $p(k_0+1)=p(k_0)$ shows
that $\bs$ has no right special factor of length $k_0$. But this implies in turn 
that $\bs$ has no right special factor of length $k_0+1$ (such a factor would
give a right special factor of length $k_0$ by removing its first letter).
Iterating shows that $\bs$ has no right special factor of length $k$, for
any $k \geq k_0$. This implies that $\bs$ is ultimately periodic ($\bs$ can be
written as a concatenation of words of length $k_0$ and each of these words must 
always be followed by the same word). 
\endpf

\bigskip

We see from Lemma~\ref{easy} above that an aperiodic sequence taking exactly 
$a$ distinct values must satisfy $p(k) \geq k + a -1$. The ``simplest'' aperiodic
sequences would thus be sequences with the smallest $p(k)$, i.e., sequences
(if any) satisfying $p(k) = k + 1$ for all $k \geq 1$. Such sequences do exist; they
are called {\em Sturmian sequences}. They are characterized in Theorem~\ref{stu} below
(see, e.g., \cite{mL02alge}). Note that Sturmian sequences admit several equivalent 
definitions and have numerous characterizations; for instance, they can be characterized 
by their palindrome or return word structure \cite{xDgP99pali, jJlV00retu}. 

\begin{theorem}\label{stu}
For any infinite word $\bs$ over $\{a,b\}$, the following properties are equivalent.
If $\bs$ satisfies these properties, then $\bs$ is called {\em Sturmian}.

\begin{itemize} 

\item The number of factors of $\bs$ of length $n$ is equal to $n+1$, for all $n \geq 1$.

\item There exist an {\em irrational} real number $\alpha > 0$ and a real number 
$\rho$, respectively called the {\em slope} and the {\em intercept} of $\bs$, 
such that $\bs$ is equal to one of the following two infinite words:
\[
  \bs_{\alpha,\rho}, ~\bs_{\alpha,\rho}^{\prime}: \NN \rightarrow \{a,b\} 
\]
defined by
\[
 \begin{matrix}
  &\bs_{\alpha,\rho}(n) = \begin{cases}
                        a    &\mbox{if} ~\lfloor(n+1)\alpha + \rho\rfloor -
                        \lfloor n\alpha + \rho\rfloor = \lfloor \alpha \rfloor \\
                        b    &\mbox{if} ~\lfloor(n+1)\alpha + \rho\rfloor -
                        \lfloor n\alpha + \rho\rfloor \neq \lfloor \alpha \rfloor
                       \end{cases} \\
  &\qquad \\
  &\bs_{\alpha,\rho}^\prime(n) = \begin{cases}
                        a    &\mbox{if} ~\lceil(n+1)\alpha + \rho\rceil -
                        \lceil n\alpha + \rho\rceil = \lfloor \alpha \rfloor \\
                        b    &\mbox{if} ~\lceil(n+1)\alpha + \rho\rceil -
                        \lceil n\alpha + \rho\rceil \neq \lfloor \alpha \rfloor
                       \end{cases}
 \end{matrix} 
\]
for $n \geq 0$ (where $\lfloor x \rfloor$ denotes the greatest integer $\leq x$ and $\lceil x \rceil$ denotes 
the least integer $\geq x$). Moreover, $\bs$ is said to be {\em characteristic Sturmian} if $\rho = \alpha$,
in which case $\bs = \bs_{\alpha, \alpha} = \bs_{\alpha, \alpha}^{\prime}$.
\end{itemize}
\end{theorem}

\begin{example}
Taking $a = 0$, $b = 1$, and $\alpha = \rho = (3 - \sqrt{5})/2$, we get the characteristic
Sturmian sequence $ 0 1 0 0 1 0 1 0 \ldots$, which is called the {\em (binary) Fibonacci 
sequence}.
\end{example}

\begin{remark}\label{slope}  
By definition it is clear that any Sturmian sequence is over a $2$-letter alphabet. 
It also follows from Lemma~\ref{easy} that Sturmian sequences are aperiodic. Note
that if we choose $\alpha$ to be rational in the above definition, we obtain (purely) 
periodic sequences, referred to as {\em periodic balanced\,} sequences -- see below. 
(Some authors also use the name {\em periodic Sturmian sequences}.) 
We will call {\em characteristic periodic balanced sequences\,} those
obtained with a {\em rational} slope $\alpha > 0$ and intercept $\rho = \alpha$ in
Theorem~\ref{stu}. Also note that the names ``slope'' and ``intercept'' refer to 
the geometric realization of Sturmian words as approximations 
to the line $y = \alpha x + \rho$ (called {\em mechanical words} in 
\cite[Chapter 2]{mL02alge}).
\end{remark}

All Sturmian sequences are ``balanced'' in the following sense.

\begin{definition} \label{D:balance} 
A finite or infinite word $w$ over $\{a,b\}$ is said to be {\em balanced} if, 
for any factors $u$, $v$ of $w$ with $|u| = |v|$, we have $||u|_{b} - |v|_{b}| \leq 1$ 
$($or equivalently $||u|_{a} - |v|_{a}| \leq 1)$. 
\end{definition}

The term ``balanced'' is relatively new; it appeared in \cite{jBpS93acha, jBpS94arem} 
(also see \cite[Chapter 2]{mL02alge}), and the notion itself dates back to 
\cite{gHmM40symb, eCgH73sequ}. In the pioneering paper \cite{gHmM40symb}, balanced 
infinite words over a $2$-letter alphabet are called ``Sturmian trajectories'' and 
belong to three classes corresponding to: Sturmian; periodic balanced; and a class 
of non-recurrent infinite words that are ultimately periodic (but not periodic), 
called {\em skew words}. 
That is, the family of balanced infinite words consists of the Sturmian words
and periodic balanced words (which are recurrent), and the (non-recurrent) skew 
infinite words, the factors of which are balanced. In particular, we have the 
following result due to Morse and Hedlund \cite{gHmM40symb}, and Coven and Hedlund 
\cite{eCgH73sequ} (see also \cite[Theorem~2.1.3]{mL02alge}):

\begin{theorem}\label{balanced}
A binary sequence is Sturmian if and only if it is balanced and aperiodic.
\end{theorem}

\begin{note}
A description of skew words is given in part (ii) of Theorem~\ref{T:fine}. Simple
examples are infinite words of the form $a^{\ell} b a^{\omega}$, where $\ell \in {\mathbb N}$.
\end{note}

It is important to note that a finite word is {\em finite Sturmian} 
(i.e., a factor of some Sturmian word) if and only if it is balanced 
\cite[Chapter 2, Proposition~2.1.17]{mL02alge}. Accordingly, 
balanced infinite words are precisely the infinite words whose factors are finite Sturmian. 
This concept is generalized in \cite{aGjJgP06char} by showing that the set of all infinite 
words whose factors are {\em finite episturmian} consists of the (recurrent) episturmian 
words and the (non-recurrent) {\em episkew} infinite words (i.e., non-recurrent infinite 
words, all of whose factors are finite episturmian), see Section~\ref{episkew}. 

\medskip

For a comprehensive introduction to Sturmian words, see for instance 
\cite{jAjS03auto, mL02alge, nP02subs} and references therein. Also see 
\cite{sBpS94orde, aHrT00char, gP05mors, rT96onco, rT98inte} for further work on skew words.

\medskip

We end this section with a simple and useful proposition which deserves to 
be better known. Its two parts were suggested several years ago to JPA 
in the case of binary sequences by J. Cassaigne and J. Berstel respectively 
(private communications).

\begin{proposition}\label{suppl}
Let $\bs$ be a sequence taking exactly $a \geq 2$ distinct values and let $p(k)$ 
be the number of distinct factors of length $k$ of $\bs$.
\begin{itemize}
\item[(i)]  If $\bs$ is aperiodic and admits at most one left special factor of
            each length, then one has $k + a - 1 \leq p(k) \leq (a - 1)k +1$ for all
            $k \geq 1$. In particular an aperiodic binary sequence which has at
            most one left special factor of each length is Sturmian.
\item[(ii)] If there exists $k_0 \geq 1$ such that $p(k) = k + a - 1$ for all
            $k \geq k_0$, then $p(k) = k + a - 1$ for all $k \geq 1$. In particular
            if a binary sequence satisfies $p(k) = k + 1$ for all $k$ larger than
            some $k_0$, then it is Sturmian.
\end{itemize}
\end{proposition}

\proof

(i).
     Using part (iv) of Lemma~\ref{easy}, we have $p(k) \geq k + a - 1$ for all $k \geq 1$, 
     since $\bs$ is aperiodic. On the other hand, erasing the first letter of all factors 
     of $\bs$ of length $k+1$ gives all factors of length $k$. There is at most one of 
     these factors of length $k$ which can be obtained from distinct factors of length
     $k+1$ (since $\bs$ admits at most one left special factor of length $k$), and if 
     so there can be at most $a$ such distinct factors of length $k+1$ (since a left 
     special factor can be extended on the left by at most $a$ letters).  
     Hence $p(k+1) - p(k) \leq a-1$ for all $k \geq 1$. By telescopic summation,
     this implies $p(k) \leq (a-1)(k-1) + p(1) = (a-1)(k-1) + a = ak - k + 1$.

(ii).
     Let $k_1$ be the least integer $\geq 1$ such that for all $k \geq k_1$, one
     has  $p(k) = k + a - 1$. Suppose that $k_1 > 1$, and let $\ell := k_1 -1$.
     Then $p(\ell) \neq \ell + a - 1$.
     But $p(\ell) \leq p(k_1) = k_1 + a - 1 = \ell + a$. Hence either
     $p(\ell) = \ell + a$, or $p(\ell) \leq \ell + a - 2$. In either case $\bs$ 
     would be ultimately periodic (by Lemma~\ref{easy} (i), resp.\ by part (iv) of
     Lemma~\ref{easy}), a contradiction. Hence $k_1 = 1$ and the claim 
     about Sturmianicity follows from Theorem~\ref{stu}.
\endpf

\subsection{Episturmian sequences} \label{SS:episturmian}

It is well known that the set of factors of any Sturmian sequence is closed
under reversal, i.e., if $u$ is a factor of a Sturmian sequence $\bs$, then
its reversal $\tilde u$ is also a factor of $\bs$ (e.g., see \cite{fM89infi} or
\cite[Proposition 2.1.19]{mL02alge}). In fact:

\begin{theorem}
An aperiodic binary sequence $\bs$ is Sturmian if and only if $F(\bs)$ is
closed under reversal and $\bs$ admits exactly one left special factor of
each length.
\end{theorem}

\begin{proof} Let $\bs$ be an aperiodic binary sequence. First suppose that
$\bs$ is Sturmian. For a proof of the fact that $F(\bs)$ is closed under
reversal, see \cite{fM89infi} or \cite[Proposition 2.1.19]{mL02alge}. Now we
will show that $\bs$ has {\em exactly} one left special factor of each length.

Let $p(n)$ denote the number of factors of length $n$ of $\bs$. Since $F(\bs)$ 
is closed under reversal, a factor of $\bs$ is left special (resp.~right special) 
if and only if its reversal is right special (resp.~left special). Hence, for all 
$n \geq 1$, the difference $p(n+1) - p(n)$ is equal to the number of left special 
factors of $\bs$ of length $n$. 
Therefore, since $p(n+1) - p(n) = 1$ for all $n \geq 1$ (by Theorem~\ref{stu}),
$\bs$ admits exactly one left special factor (or equivalently, right special 
factor) of each length.

The converse follows immediately from part (i) of Proposition~\ref{suppl}.
\endpf
\end{proof}

Inspired by results of this flavour, Droubay, Justin, and Pirillo
\cite{xDjJgP01epis, jJgP02epis} introduced the following natural
generalization of Sturmian sequences on an arbitrary finite alphabet $\cA$.

\begin{definition} {\em \cite{xDjJgP01epis}} An infinite word $\bt \in
\cA^\omega$ is said to be {\em episturmian} if its set of factors $F(\bt)$
is closed under reversal and $\bt$ admits at most one left special factor
(or equivalently, right special factor) of each length.
\end{definition}

\begin{note}
When $\cA$ is a 2-letter alphabet, this definition gives the Sturmian words
as well as the periodic balanced words.
\end{note}

In the seminal paper \cite{xDjJgP01epis}, episturmian words were defined
as an extension of {\em standard episturmian words}, which were first
introduced as a generalization of characteristic Sturmian words using {\em
iterated palindromic closure} (a construction due to de~Luca
\cite{aD97stur}).

The {\em palindromic right-closure} $w^{(+)}$ of a finite word $w$ is the
(unique) shortest palindrome beginning with $w$ (see \cite{aD97stur}). 
More precisely, if $w = uv$ where $v$ is the longest palindromic suffix of $w$, 
then $w^{(+)} = uv\rev{u}$. For example, {\tt (tie)}$^{(+)} =$ {\tt tie\thinspace it}. 
The {\em iterated palindromic closure function} \cite{jJ05epis}, denoted by $Pal$, 
is defined recursively as follows. Set $Pal(\empt) = \empt$ and, for any word
$w$ and letter $x$, define $Pal(wx) := (Pal(w)x)^{(+)}$. For instance,
$Pal(abc) = (Pal(ab)c)^{(+)} = (abac)^{(+)} = abacaba$. Note that $Pal$ is 
injective; and moreover, it is clear from the definition that $Pal(w)$ is a prefix of 
$Pal(wx)$ for any word $w$ and letter $x$. Hence, if $v$ is a prefix of $w$, then 
$Pal(v)$ is a prefix of $Pal(w)$.

\newpage
\begin{theorem} \label{T:xDjJgP01epis} {\em \cite{xDjJgP01epis}}
For an infinite word $\bs \in \cA^\omega$, the following properties are
equivalent.
\begin{itemize}
\item[(i)] There exists an infinite word $\Delta = x_1x_2x_3\ldots$ ($x_i
\in \cA$), called the {\em directive word} of $\bs$, such that $\bs =
\lim_{n\rightarrow\infty}Pal(x_1x_2\cdots x_n)$.
\item[(ii)] $F(\bs)$ is closed under reversal and all of the left special
factors of $\bs$ are prefixes of it.
\end{itemize}
An infinite word $\bs$ satisfying the above properties is said to be 
{\em standard episturmian} (or {\em epistandard} for short).
\end{theorem}

The above characterization of epistandard words extends to the case of an
arbitrary finite alphabet a construction given in \cite{aD97stur} for all
characteristic Sturmian words.

\begin{example} \label{Ex:Trib}
The epistandard word $\br$ directed by $\Delta = (abc)^\omega$ is known as
the {\em Tribonacci word}; it begins in the following way:
\[
  \br =
\underline{a}\underline{b}a\underline{c}aba\underline{a}bacaba\underline{b}
acabaabacaba\underline{c}abaabaca
\cdots~,
\]
where each palindromic prefix $Pal(x_1\cdots x_{n-1})$ is followed by an
underlined letter $x_{n}$. More generally, for $k\geq 2$, the {\em
$k$-bonacci word} is the epistandard word over $\{a_1, a_2, \ldots, a_k\}$
directed by $(a_1a_2\cdots a_k)^\omega$.
\end{example}

\begin{remark}
In \cite{xDjJgP01epis}, Droubay {\it et al.} proved that an infinite word
$\bt$ is episturmian if and only if $F(\bt) = F(\bs)$ for some epistandard
word $\bs$. They also proved that episturmian words are uniformly recurrent;
hence any such infinite word is either (purely) periodic or aperiodic. The
aperiodic episturmian words are precisely the episturmian words that admit
exactly one left special factor of each length. In fact, an epistandard word
$\bs$ (and hence any episturmian word with the same set of factors $\bs$) is
periodic if and only if exactly one letter occurs infinitely often in the
directive word of $\bs$ (see \cite[Proposition 2.9]{jJgP02epis}).
\end{remark}

The notion of a {\em directive word} (as defined for epistandard words in
Theorem~\ref{T:xDjJgP01epis}) extends to {\em all} episturmian words with
respect to {\em episturmian morphisms}, which play a central role in the
study of these words. Introduced first as a generalization of {\em Sturmian
morphisms}, Justin and Pirillo \cite{jJgP02epis} showed that episturmian
morphisms are exactly the morphisms that preserve the aperiodic episturmian
words (i.e., the morphisms that map aperiodic episturmian words onto
aperiodic episturmian words). Such morphisms naturally generalize to any
finite alphabet the {\em Sturmian morphisms} on two letters. A morphism
$\varphi$ is said to be {\em Sturmian} if $\varphi(\bs)$ is Sturmian for any
Sturmian word $\bs$. The set of Sturmian morphisms over $\{a,b\}$ is closed
under composition, and consequently it is a submonoid of the endomorphisms
of $\{a,b\}^*$. Moreover, it is well known that the monoid of Sturmian
morphisms is generated by the three morphisms: $(a\mapsto ab, b\mapsto a)$,
$(a\mapsto ba. b \mapsto a)$, $(a\mapsto b, b \mapsto a)$ and that Sturmian
morphisms are precisely the morphisms that map Sturmian words onto Sturmian
words (see \cite{jBpS93acha, fMpS93morp}; also see Section~\ref{stumorp} later).

By definition (see \cite{xDjJgP01epis, jJgP02epis}), the monoid of all 
{\em episturmian morphisms} is generated, under composition, by all the
morphisms:
\begin{itemize}
\item $\psi_a$:  $\psi_a(a) = a$, $\psi_a(x) = ax$ for any letter $x \neq a$;
\item $\bar \psi_a$:  $\bar \psi_a(a) = a$, $\bar \psi_a(x) = xa$ for any
letter $x \neq a$;
\item $\theta_{ab}$: exchange of letters $a$ and $b$.
\end{itemize}

Moreover, the monoid  of so-called {\em epistandard morphisms} is generated
by all the $\psi_a$ and the $\theta_{ab}$, and the monoid of {\em pure
episturmian morphisms} (resp.~{\em pure epistandard morphisms}) is generated
by the $\psi_a$ and $\bar{\psi}_a$ only (resp.~the $\psi_a$ only). The
monoid  of the {\em permutation morphisms} (i.e., the morphisms $\varphi$
such that $\varphi(\cA) = \cA$) is generated by all the $\theta_{ab}$.

As shown in \cite{jJgP02epis}, any episturmian word is the image of another
episturmian word by some pure episturmian morphism and any episturmian word
can be infinitely decomposed over the set of pure episturmian morphisms.
This last property allows an episturmian word to be defined by one of its
morphic decompositions or, equivalently, by a certain {\em spinned directive
word}, which is an infinite sequence of rules for decomposing the given
episturmian word by morphisms. See \cite{aGfLgR08dire, jJgP04epis} for
recent work concerning directive words of episturmian words.

\begin{remark} \label{R:episturmian-orbit}
The {\em shift-orbit} of an infinite word $\bx \in \cA^\omega$ is the set 
$\Orbit(\bx) = \{T^i(\bx), \, i \geq 0\}$ and its {\em closure} is given by 
$\Closure(\bx) = 
\left\{\by \in \cAw, \, \Pref(\by) \subseteq \bigcup_{i\geq0}\Pref(T^i(\bx))\right\}$, 
where $\Pref(w)$ denotes the set of prefixes of a finite or infinite word $w$. 
Note that for any infinite word $\bt$ and $\bx \in \Closure(\bt)$,
$F(\bx) \subseteq F(\bt)$. If, moreover, $\bt$ is uniformly recurrent,
then it follows that for each $n \geq 1$, $F_n(\bx) = F_n(\bt)$, and
hence $F(\bx) = F(\bt)$ for any $\bx \in \Closure(\bt)$ (see for
instance \cite[Proposition~5.1.10]{nP02subs} or
\cite[Proposition~1.5.9]{mL02alge}). This implies that $\Closure(\bx)
= \Closure(\bt)$ for any $\bx \in \Closure(\bt)$; in other words,
$\Closure(\bt)$ is a {\em minimal dynamical system} (see, e.g.,
\cite{mL02alge, nP02subs}). Accordingly, since episturmian words are
uniformly recurrent, the closure of the shift-orbit of any episturmian
$\bt$ is a minimal dynamical system; in particular, $\Closure(\bt)$
consists of all the episturmian words with the same set of factors as
$\bt$ (see, e.g., \cite{rRlZ00agen}).

Note that if $\bt$ is aperiodic, then $\Closure(\bt)$ contains a unique 
epistandard word with the same set of factors as $\bt$, whereas if $\bt$ 
is periodic, $\Closure(\bt)$ contains two different epistandard words  
(see for instance \cite{aGjJ07epis, aGfLgR08dire}).
\end{remark}

\subsubsection{Strict episturmian words}

\begin{definition} An epistandard word $\bs$ (or any episturmian word with
the same set of factors as $\bs$) is said to be {\em strict} if every letter
in the alphabet of $\bs$ occurs infinitely often in its directive word.
\end{definition}

Strict episturmian words on $k$ letters are often said to be {\em $k$-strict}; 
these words have $(k-1)n+1$ distinct factors of length $n$ for all $n \geq 1$ 
(as proven in \cite[p.~549]{xDjJgP01epis}) and they coincide with the $k$-letter 
{\em Arnoux-Rauzy sequences} introduced in \cite{pAgR91repr}. In particular, 
the $2$-strict episturmian words are exactly the Sturmian words since these 
words have $n+1$ distinct factors of length $n$ for each $n\geq 1$ (recall 
Theorem~\ref{stu}).

Note that any episturmian word takes the form $\varphi(\bt)$ with $\varphi$ an 
episturmian morphism and $\bt$ an Arnoux-Rauzy sequence (or strict episturmian 
word). In this sense, episturmian words are only a slight generalization of 
Arnoux-Rauzy sequences. For example, the family of episturmian words on three 
letters $\{a,b,c\}$ consists of the Arnoux-Rauzy sequences over $\{a,b,c\}$, 
the Sturmian words over $\{a,b\}$, $\{b,c\}$, $\{a,c\}$ and their images under 
episturmian morphisms on $\{a,b,c\}$, and periodic infinite words of the form 
$\varphi(x)^\omega$ where $\varphi$ is an episturmian morphism on $\{a,b,c\}$ 
and $x \in \{a,b,c\}$.

\subsubsection{Episkew words}\label{episkew}

A finite word $w$ is said to be {\em finite Sturmian} (resp.~{\em finite
episturmian}) if $w$ is a factor of some infinite Sturmian
(resp.~episturmian) word.

Recall from Section~\ref{SS:Sturmian} that skew words are ultimately 
periodic (but not periodic) infinite words, all of whose factors
are finite Sturmian (or equivalently, balanced). Over a 2-letter alphabet, 
skew words constitute the family of non-recurrent balanced infinite words, 
whereas the recurrent balanced infinite words consist of the Sturmian words 
and the periodic balanced words.

Inspired by Morse and Hedlund's~\cite{gHmM40symb} skew words, 
{\em episkew words} were recently defined in \cite{aGjJgP06char} as
non-recurrent infinite words, all of whose factors are finite episturmian. 
A number of equivalent definitions of such words were given in
\cite{aGjJgP06char} (also see Theorem \ref{T:fine}, to follow).

Episkew words were first alluded to (but not explicated) in \cite{aG06acha}.
Following that paper, these words showed up again in the study of
inequalities characterizing finite and infinite episturmian words in
relation to lexicographic orderings \cite{aGjJgP06char}; in fact, as
detailed in Section \ref{SS:Pirillo}, episturmian words have extremal
properties similar to those of Sturmian words.

To learn more about episturmian and episkew words, see for instance the
recent surveys \cite{jB07stur, aGjJ07epis}.

\section{Extremal words}

Suppose the alphabet $\cA$ is totally ordered by the relation $\leq$. Then we 
can totally order $\cA^+$ by the {\em lexicographic order} $\leq$, defined as 
follows. Given two non-empty finite words $u$, $v$ on $\cA$, we have $u < v$ if 
and only if either $u$ is a prefix of $v$ (with $u \ne v$) or 
$u = xau^\prime$ and $v = xbv^\prime$, for some finite words $x$, $u^\prime$,
$v^\prime$ and letters $a$, $b$ with $a < b$. This is
the usual alphabetic ordering in a dictionary, and we say that $u$
is {\em lexicographically less} than $v$. This notion 
naturally extends to infinite words, as follows. Let $\bu =
u_0u_1u_2\cdots$ and $\bv = v_0v_1v_2\cdots$, where $u_j$, $v_j
\in \cA$. We define $\bu < \bv$ if there exists an index $i\geq0$
such that $u_j = v_j$ for all $j=0,\ldots, i-1$ and $u_{i} < v_{i}$.

Let $w$ be a finite or infinite word on $\cA$, and let $k$ be a positive integer. We let
 $\min(w | k)$ (resp.~$\max(w | k)$) denote the lexicographically smallest (resp.~greatest) 
factor of $w$ of length $k$ for the given order (where $|w|\geq k$ if $w$ is finite).  

If $w$ is infinite, then it is clear that $\min(w | k)$ and $\max(w | k)$ are prefixes 
of the respective words $\min(w | k+1)$ and $\max(w | k+1)$. So we can define, by taking 
limits, the following two infinite words (see \cite{gP05ineq}):
\[
  \min(w) = \underset{k\rightarrow\infty} {\lim}\min(w | k) \quad \mbox{and} \quad 
  \max(w) = \underset{k\rightarrow\infty}{\lim}\max(w | k).
\]  
That is, to any infinite word $\bt$ we can associate two infinite words $\min(\bt)$ and 
$\max(\bt)$ such that any prefix of $\min(\bt)$ (resp.~$\max(\bt)$) is the 
{\em lexicographically} smallest (resp.~greatest)  amongst the factors of $\bt$ of the 
same length.

For a finite word $w$ on a totally ordered alphabet $\cA$, $\min(w)$ denotes 
$\min(w | k)$ where $k$ is maximal such that all $\min (w | j)$, $j= 1,2, \dots, k$, 
are prefixes of $\min (w | k)$. In the case $\cA= \{a,b\}$, $\max(w)$ is defined 
similarly (see \cite{aGjJgP06char}).

The following definition, given in \cite{aGjJgP06char}, will be useful in the next section, 
where we survey recent work concerning extremal properties of (epi)Sturmian sequences, 
particularly inequalities characterizing such words (finite and infinite).

\begin{definition} An {\em acceptable pair} for an alphabet $\cA$ is a pair $(a, <)$ 
where $a$ is a letter in $\cA$ and $<$ is a total order on $\cA$ such that $a = \min(\cA).$
\end{definition}

\section{Extremal properties}

In 2003, Pirillo \cite{gP03ineq} (also see \cite{gP05ineq}) proved that, for infinite words 
$\bs$ on a $2$-letter alphabet $\{a,b\}$ with $a<b$, the inequality 
\begin{equation} \label{eq:characteristic}
a\bs \leq \min(\bs) \leq \max(\bs) \leq b\bs
\end{equation}
characterizes the {\em characteristic Sturmian words} and 
{\em characteristic periodic balanced words}.

\begin{remark}
Characteristic periodic balanced sequences, which correspond to the ``Sturmian'' 
sequences with {\em rational} slope $\alpha > 0$ and intercept $\rho = \alpha$ 
(see Theorem~\ref{stu} and Remark~\ref{slope}) are precisely the sequences 
of the form $(Pal(v)xy)^\omega$ where $v \in \{a,b\}^*$ and $\{x,y\} =
\{a,b\}$ (see for instance \cite{jAaG09dist, jB07stur, xDjJgP01epis}). Also note that if
$\bs$ is a characteristic Sturmian sequence, then $a\bs = \min(\bs)$
and $b\bs = \min(\bs)$. On the other hand, if $\bs$ is a
characteristic periodic balanced sequence, then either:
\begin{itemize}
\item  $a\bs < \min(\bs)$ and $b\bs = \max(\bs)$ when $\bs$ takes the
form $(Pal(v)ab)^\omega$,
\item or $a\bs = \min(\bs)$ and $\max(\bs) < b\bs$ when $\bs$ takes
the form $(Pal(v)ba)^\omega$.
\end{itemize}
For example, the characteristic periodic balanced
sequence $\bs := (Pal(ab)ab)^\omega = (abaab)^\omega$ satisfies
\[
a\bs = a(abaab)^\omega < \min(\bs) = (aabab)^\omega \quad \mbox{and}
\quad b\bs = b(abaab)^\omega = \max(\bs),
\]
whereas $\bs' := (Pal(ab)ba)^\omega = (ababa)^\omega$ satisfies
\[
a\bs' = a(ababa)^\omega = \min(\bs') \quad \mbox{and} \quad 
\max(\bs') = (babaa)^\omega < b\bs' = b(ababa)^\omega.
\]
More generally, given two characteristic periodic balanced sequences
$\bs$, $\bs'$ of the form $\bs = (Pal(v)ab)^\omega$ and $\bs' =
(Pal(v)ba)^\omega$ for some $v \in \{a,b\}^*$, we have
\[
  \min(\bs) = \min(\bs') = (aPal(v)b)^\omega \quad \mbox{and} \quad
\max(\bs) = \max(\bs') = (bPal(v)a)^\omega.
 \]
 See \cite{jAaG09dist, gP05ineq} for more details.
\end{remark}

\medskip

The preceding result of Pirillo concerning characteristic Sturmian words and 
characteristic periodic balanced words (property~\eqref{eq:characteristic}) 
encompasses Theorem \ref{T:intro} -- one of the key properties underlying the 
main theorem in Bugeaud and Dubickas' paper \cite{yBaD05frac}. In fact, as 
mentioned previously, Theorem~\ref{T:intro} was known much earlier -- in 1993, 
Berstel and S{\ee}bold~\cite{jBpS93acha} (as well as Borel and Laubie 
\cite{jBfL93quel}) proved one direction of the theorem, namely that 
characteristic Sturmian words satisfy \eqref{eq:characteristic}. This Sturmian 
extremal property also resurfaced in 2001, under a different guise, in a paper 
of S.~Gan \cite{sG01stur}. However, it seems that P. Veerman \cite{pV87symb} 
was actually the first to prove \eqref{eq:characteristic} for Sturmian sequences 
in 1987, albeit from a symbolic dynamical perspective and in an implicit way. 
A year prior, Veerman had already proved that characteristic Sturmian sequences 
have the above extremal property \cite[Theorem 2]{pV86symb}; it was not until 
\cite[Theorem 2.1]{pV87symb} that he proved the equivalence. Motivated by the 
combinatorics of the Mandelbrot set, Bullett and Sentenac \cite{sBpS94orde} 
reproved these results of Veerman, in the language of ordered sets.

In this section, we shall first discuss the combinatorial work of Pirillo and others in 
relation to the inequalities \eqref{eq:characteristic} and their generalizations. Following 
this, we will consider in more detail the earlier work by Berstel and S{\ee}bold 
\cite{jBpS93acha}, Gan \cite{sG01stur}, and Veerman \cite{pV86symb, pV87symb}. 

\subsection{Pirillo's work continued} \label{SS:Pirillo}

Continuing his work in relation to the inequalities \eqref{eq:characteristic}, Pirillo 
\cite{gP05ineq} proved further that, in the case of an arbitrary finite alphabet $\cA$, 
an infinite word $\bs$ on $\cA$ is {\em epistandard} if and only if, for any 
acceptable pair $(a,<)$, we have 
\begin{equation} \label{eq:gP05ineq}
a\bs \leq \min(\bs). 
\end{equation} 
Moreover, $\bs$ is a strict epistandard word if and only if \eqref{eq:gP05ineq} 
holds with strict equality for any order \cite{jJgP02onac}.

In a similar spirit, Pirillo \cite{gP05mors} defined {\em fine words} over two letters; 
that is, an infinite word $\bt$ over a $2$-letter alphabet $\{a,b\}$ ($a < b$) 
is said to be {\em fine} if  $(\min(\bt), \max(\bt)) = (a\bs, b\bs)$ for some infinite 
word $\bs$. These infinite words were characterized in \cite{gP05mors} by showing that fine 
words on $\{a,b\}$ are exactly the {\em Sturmian} and {\em skew} infinite words
(see Section~\ref{SS:Sturmian}). Specifically:

\begin{theorem} \label{T:gP05mors} Let $\bt$ be an infinite word over $\{a,b\}$. 
The following properties are equivalent:  
\begin{itemize}
\item[(i)] $\bt$ is fine,
\item[(ii)] either $\bt$ is Sturmian, or $\bt$ is an ultimately periodic 
                   (but not periodic) shift of an infinite word of the form
                   $\mu(x^{\ell} y x^{\omega})$ for some $\ell \in {\mathbb N}$, 
                   where $\mu$ is a pure standard morphism on $\{a, b\}$ and
                   $\{x, y\} = \{a, b\}$ (these words are the skew words).
\end{itemize}  
\end{theorem}
In other words, a fine word over two letters is either a Sturmian word or an 
ultimately periodic (but not periodic) infinite word, all of whose factors are Sturmian. 

Pirillo \cite{gP05mors} remarked that perhaps his characterization of fine words 
could be generalized to an arbitrary finite alphabet; indeed, Glen \cite{aG06acha} 
soon generalized this result by extending Pirillo's definition of fine words to more 
than two letters. That is:

\begin{definition} \label{D:fine} {\rm \cite{aG06acha}}
An infinite word $\bt$ on $\cA$ is said to be {\em fine} if  there exists an infinite word 
$\bs$ such that $\min(\bt) = a\bs$ for any acceptable pair $(a,<)$.
\end{definition} 

\begin{note}
It is easy to see that Pirillo's original $2$-letter definition of a fine word is a special 
instance of the above definition. Certainly, as there are only two lexicographic orders on 
words over a $2$-letter alphabet, it follows from Definition \ref{D:fine} that a fine word 
$\bt$ over $\{a,b\}$ ($a< b$) satisfies $(\min(\bt), \max(\bt)) = (a\bs, b\bs)$ for some 
infinite word $\bs$.
\end{note}

Glen \cite{aG06acha} characterized these generalized fine words (given in
Definition \ref{D:fine}) by showing that such an infinite word is either 
a {\em strict} episturmian word or a {\em strict episkew word}. More precisely:

\newpage
\begin{theorem} \label{T:fine} {\rm \cite{aG06acha}} Let $\bt$ be an infinite word with 
$\mbox{{\em Alph}}(\bt) = \cA$. Then, $\bt$ is fine if and only if one of the 
following holds:
\begin{itemize}
\item[(i)]  $\bt$ is an $\cA$-strict episturmian word; 
\item[(ii)] $\bt$ is non-recurrent and takes the form $\mu(x\bs)$ where $x$
            is a letter, $\bs$ is a strict epistandard word on
            $\cA\setminus\{x\}$, and $\mu$ is a pure episturmian morphism on $\cA$.
\end{itemize} 
\end{theorem}

\begin{remark}
Note that part (ii) of Theorem~\ref{T:fine} gives the form of
so-called {\em strict episkew words}; it is slightly simpler to what
was originally given in \cite{aGjJgP06char}, thanks to
Richomme (private communication). Also note that strict episkew words on a
$2$-letter alphabet are precisely the skew words (see \cite{aGjJ07epis}).
One can also compare Theorem~\ref{T:fine} with Theorem~\ref{T:gP05mors}.
A simple example of an episkew word is 
$c \mathbf{f} := c a b a a b a b a a b a \ldots$, where $\mathbf{f}$ is the 
Fibonacci sequence on $\{a, b\}$.
\end{remark}

\begin{example} \cite{aG06acha} Let $\cA = \{a,b,c\}$ with $a<b<c$. Let $\bbf$ denote 
the infinite Fibonacci word over $\{a,b\}$, i.e., the epistandard
word directed by $(ab)^\omega$. Then, the following infinite words are fine. 
\begin{itemize}
\item $\bbf = abaababaabaaba\cdots$ \smallskip
\item $c\bbf = \underline{c}abaababaabaaba\cdots$ \smallskip
\item $\rev{\bbf_4}c\bbf = aaba\underline{c}abaababaabaaba\cdots$ \smallskip
\item $\psi_a(\bbf) = aabaaabaabaaabaaaba\cdots$ \smallskip
\item $\psi_c(c\bbf) = \underline{c}cacbcacacbcacbcacacbcacacbca\cdots$ \smallskip 
\item $\psi_c(\rev{\bbf_4}c\bbf) = cacacbca\underline{c}cacbcacacbcacbcacacbcaca\cdots$ 
\smallskip
\end{itemize}
Let us note, for example, that $\psi_c(\bbf)$ is {\bf not} fine since it is a 
{\em non-strict} epistandard word. That is, $\psi_c(\bbf)$ is an epistandard word with 
directive word $c(ab)^\omega$, so it is not strict, nor does it take the second form given 
in Theorem \ref{T:fine}. 
\end{example}
 
Continuing this work, Glen, Justin, and Pirillo \cite{aGjJgP06char} recently proved 
new characterizations of {\em finite} Sturmian and episturmian words via lexicographic 
orderings. As a consequence, they were able to characterize by lexicographic order all 
episturmian words in a {\em wide sense} (episturmian and episkew infinite words). 
Similarly, they characterized by lexicographic order all balanced infinite words on a 
$2$-letter alphabet; in other words, all Sturmian, periodic balanced, and skew infinite 
words, the factors of which are (finite) Sturmian.

In the finite case:

\begin{theorem}\label{T:aGjJgP06char} {\rm \cite{aGjJgP06char}} 
A finite word $w$ on $\cA$ is episturmian if and only if there exists a finite 
word $u$ on $\cA$ such that, for any acceptable pair  $(a, <)$, we have 
\begin{equation}
au_{|m|-1} \le m \label{e2} 
\end{equation} 
where $m= \min(w)$ for the considered order. \qed
\end{theorem}

A corollary of Theorem \ref{T:aGjJgP06char} is the following new characterization of 
finite Sturmian words (i.e., finite balanced words).

\begin{corollary}\label{Cor:aGjJgP06char-finite} {\rm \cite{aGjJgP06char}}  
A finite word $w$ on $\cA = \{a,b\}$, $a<b$, is not Sturmian (in other words, not balanced)  
if and only if there exists a finite word $u \in \{a,b\}^*$ such that $aua$ is a prefix of 
$\min(w)$ and $bub$ is a prefix of $\max(w)$. 
\qed
\end{corollary}

In the infinite case, a characterization of episturmian words in the  {\em wide sense} 
follows almost immediately from Theorem \ref{T:aGjJgP06char}. That is:

\begin{corollary}\label{Cor:aGjJgP06char-infinite} {\rm \cite{aGjJgP06char}}  
An infinite word $\bt$ on $\cA$ is episturmian in the wide sense (i.e., episturmian 
or episkew) if and only if there exists an infinite word $\bu$ on $\cA$ such that
\begin{equation*}a\bu \le \min(\bt) 
\end{equation*}
for any acceptable pair  $(a, <)$.
\end{corollary}

Consequently, an infinite word $\bs$ on $\{a,b\}$ ($a< b$) is balanced (i.e., Sturmian,
periodic balanced, or skew) if and only if there exists an infinite word $\bu$ on
$\{a, b\}$ such that 
\begin{equation}
a\bu \le \min(\bs) \le \max(\bs) \le b\bu. \label{e1}
\end{equation}

For any sequence $\bs$, $\max(\bs)$ is the same as $\sup\{T^k(\bs), k\geq 0\}$, 
and similarly $\min(\bs) = \inf\{T^k(\bs), k\geq 0\}$, where the infimum and supremum are taken with respect to the lexicographic order. The preceding
result therefore shows that a sequence $\bs$ in $\{0,1\}^\omega$ is balanced if
and only if there exists a sequence $\bu \in \{0,1\}^\omega$ such that $0\bu \leq
T^k(\bs) \leq 1\bu$  for all $k \geq 0$. In particular, a sequence $\bs$ on
$\{0,1\}$ being Sturmian is equivalent to $\bs$ being aperiodic and the
existence of a sequence $\bu$ on $\{0,1\}$ such that $0\bu \leq T^k(\bs)
\leq 1\bu$. Moreover, it follows from the proof of Theorem~\ref{T:gP05mors}
(or Theorem~\ref{T:fine}) that $\bu$ is the unique characteristic Sturmian
sequence having the same slope as $\bs$. This is exactly Theorem 2.1. For
the sake of completeness, we give a direct proof below.

\bigskip

\noindent{\it Direct proof of Theorem~$\ref{P:JPA}$.\ \ }
Let $\bs$ be an aperiodic sequence on $\{0, 1\}$. First suppose 
that $\bs$ is a Sturmian sequence. Since it contains both
$0$'s and $1$'s, there exist two binary sequences $\bx$ and $\by$
such that $0 \bx := \inf\{T^k(\bs), \ k \geq 0\}$ and 
$1 \by := \sup\{T^k(\bs), \ k \geq 0\}$. We claim that $\bx \geq \by$.
Namely, if $\bx < \by$, there exist a (possibly empty) word $w$ and 
two infinite sequences $\bx'$ and $\by'$ such that $\bx = w 0 \bx'$ 
and $\by = w 1 \by'$. Hence $0\bx = 0 w 0 \bx'$ and $1\by = 1 w 1 \by'$.
Since any factor of $\inf\{T^k(\bs), \ k \geq 0\}$ and of 
$\sup\{T^k(\bs), \ k \geq 0\}$ is a factor of $\bs$, we have that both 
$0w0$ and $1w1$ are factors of $\bs$. Hence $\bs$ is unbalanced (see
Definition~\ref{D:balance} and the comments following it), but is was
supposed Sturmian, a contradiction (Theorem~\ref{balanced}). 
Thus $\bx \geq \by$, and hence
$$
\forall k \geq 0, \ 0 \bx \leq T^k(\bs) \leq 1 \by \leq 1 \bx.
$$

\medskip

Now suppose that $\bs$ has the property that there exists a binary sequence 
$\bu$ such that
\begin{equation}\label{ineq}
\forall k \geq 0, \ 0 \bu \leq T^k(\bs) \leq 1 \bu.
\end{equation}
Let $z$ be a left special factor (if any) of $\bs$, and let $z'$ be the prefix of 
$\bu$ that has the same length as $z$. Since $0z$ and $1z$ are both factors of $\bs$, 
there exist two integers $\ell_1$ and $\ell_2$ such that $T^{\ell_1}(\bs)$ begins with 
$0z$ and $T^{\ell_2}(\bs)$ begins with $1z$. We deduce from the inequalities~(\ref{ineq}) 
with $k = \ell_1$ (resp.\ $\ell_2$) that
$$
0 z' \leq 0z \ \ \mbox{\rm and} \ \ 1 z \leq 1 z'.
$$
This implies
$$
z' \leq z \ \ \mbox{\rm and} \ \ z \leq z'
$$
hence $z = z'$. Thus $\bs$ has at most one left special factor of each length.
Hence $\bs$ is Sturmian (Proposition~\ref{suppl}), and its left special factors 
are exactly the prefixes of $\bu$.

\medskip

This implies furthermore that $\bu$ belongs to the closure of the shift-orbit of $\bs$,
hence it is Sturmian. But the prefixes of $0 \bu$ and $1 \bu$ are also 
factors of $\bs$. Hence $0 \bu$ and $1 \bu$ are also in the closure of the shift-orbit of $\bs$,
thus Sturmian. This implies that $\bu$ is Sturmian characteristic (see, e.g., 
\cite[Proposition~2.1.22]{mL02alge}). Thus $\bu$ is the (unique) characteristic 
Sturmian sequence having the same slope as $\bs$. \endpf

\begin{remark}\label{2.1-2.2}
We noted in the Introduction that Theorem~\ref{T:intro} can be easily deduced from
Theorem~\ref{P:JPA}. Actually Theorem~\ref{P:JPA} can also be deduced from 
Theorem~\ref{T:intro}: it suffices to remember that the closure of the shift-orbit 
of a characteristic Sturmian sequence $\bu$ is exactly the set of all Sturmian 
sequences having the same slope as $\bu$ (see for instance 
\cite[Proposition~2.1.25]{mL02alge}), and all of these Sturmian
sequences have the same set of factors (\cite[Proposition~2.1.18]{mL02alge}, or 
\cite{fM89infi}). See also Remark~\ref{R:orbit-closure} later. 
\end{remark}

Recently, Richomme \cite{gR07aloc} proved that episturmian words can be characterized 
via a nice ``local balance property''. That is: 

\begin{theorem} \label{T:gR07aloc} {\em \cite{gR07aloc}} For a recurrent infinite word 
$\bt \in \cAw$, the following assertions are equivalent:
\begin{enumerate}
\item[(i)] $\bt$ is episturmian; 
\item[(ii)] for each factor $u$ of $\bt$, there exists a letter $a$ such that 
$\cA u\cA\cap F(\bt) \subseteq au\cA\cup\cA ua$;
\item[(iii)] for each palindromic factor $u$ of $\bt$, there exists a letter $a$ 
such that $\cA u\cA\cap F(\bt) \subseteq  au\cA \cup \cA ua$. 
\end{enumerate}
\end{theorem} 

Roughly speaking, the above theorem says that for any factor $u$ of a given episturmian 
word $\bt$, there exists a unique letter $a$ such that  every occurrence of $u$ in $\bt$ 
is immediately preceded or followed by $a$ in $\bt$. When $|\cA| = 2$, property (ii) 
of Theorem~\ref{T:gR07aloc} is equivalent to the definition of balance. Indeed, 
Coven and Hedlund \cite{eCgH73sequ} stated that an infinite word $\bs$ over $\{a, b\}$ is 
not balanced if and only if there exists a palindrome $u$ such that $aua$ and $bub$ are 
both factors of $\bs$. As pointed out in \cite{gR07aloc}, this property can be rephrased 
as follows: an infinite word $\bs$ is Sturmian if and only if $\bs$ is aperiodic and, for 
any factor $u$ of $\bs$, the set of factors belonging to $\cA u \cA$ is a subset of 
$au\cA \cup \cA u a$ or a subset of $bu\cA \cup \cA ub$.

\begin{remark} Recall that the set of all infinite words in $\cA^\omega$ having episturmian 
factors consists of the (recurrent) episturmian words and the (non-recurrent) episkew words 
in $\cA^\omega$. Therefore, since properties (ii) and (iii) in Theorem~\ref{T:gR07aloc} 
concern only factors, one readily deduces that these properties in fact characterize the 
episturmian and episkew words in $\cA^\omega$. So the hypothesis of recurrence in the 
statement of the theorem restricts attention to episturmian words only.
\end{remark}

We will now use Theorem~\ref{T:gR07aloc} to give an alternative (simpler) proof the following 
analogue of Theorem~2.1 for episturmian sequences, which was originally proved in \cite{aG07orde} 
(also see~\cite{aGjJgP06char}). This result, in particular, gives a more precise version of 
Corollary~\ref{Cor:aGjJgP06char-infinite} under the hypothesis of recurrence.

\begin{theorem}\label{AG}
 A recurrent infinite word $\bt$ on $\cA$ is episturmian if and only if 
there exists an infinite word $\bu$ on $\cA$ such that, for any acceptable pair $(a, <)$,
\[
a\bu \leq T^i(\bt)\quad \mbox{for all $i \geq 0$}.
\]
Moreover, if $\bt$ is aperiodic, then $\bu$ is the unique epistandard word with the same set 
of factors as $\bt$ (i.e., the unique epistandard word in the closure of the shift-orbit of 
$\bt$),  and for any acceptable pair $(a,<)$, $a\bu = \inf\{T^k(\bt),\ k \geq 0\}$ if and only 
if the letter $a$ occurs infinitely often in the directive word of $\bu$.
\end{theorem}

\begin{proof}
Let $\bt$ be a recurrent infinite word on $\cA$. 

First suppose that $\bt$ is episturmian. Let $x$ be a letter in $\cA$ and consider two different 
total orders $<_1$ and $<_2$ on $\cA$ such that $(x,<_1)$ and $(x,<_2)$ are acceptable pairs. 
Then there exist infinite words $\bu$ and $\bv$ on $\cA$ such that 
\begin{equation} \label{eq:inf-1}
  x\bu = \mbox{$\inf_1$}\{T^k(\bt), k\geq 0\} \quad \mbox{for the total order $<_1$ on $\cA$},
\end{equation}
and
\begin{equation} \label{eq:inf-2}
  x\bv = \mbox{$\inf_2$}\{T^k(\bt), k\geq 0\} \quad \mbox{for the total order $<_2$ on $\cA$}.
\end{equation}
(Here, $\inf_i$ denotes the infimum with respect to the order $<_i$ for $i = 1,2$.)  
We will show that $\bu = \bv$. By equations~\eqref{eq:inf-1} and \eqref{eq:inf-2}, we have
\[
x\bu \leq_1 x\bv \quad \mbox{and} \quad x\bv \leq_2 x\bu.
\]
Hence, if $u$ and $v$ are prefixes of the respective words $\bu$ and $\bv$ with $|u| = |v|$, 
then we have $u \leq_1 v$ and $v \leq_2 u$. This implies that $u = v$, and therefore 
$\bu = \bv$. Hence, for a given letter $x$ in $\cA$, there exists a unique infinite word 
$\bu$ on $\cA$ such that
\begin{equation} \label{eq:inf-1*}
x\bu = \mbox{$\inf_x$}\{T^k(\bt), k\geq 0\} \quad \mbox{for any acceptable pair $(x,<_x)$}. 
\end{equation}
Now consider another letter $y$ in $\cA\setminus\{x\}$. By what precedes, we know there 
exists a unique infinite word $\bv$ on $\cA$ such that
\begin{equation} \label{eq:inf-2*}
y\bv =  \mbox{$\inf_y$}\{T^k(\bt), k\geq 0\} \quad \mbox{for any acceptable pair $(y,<_y)$}.
\end{equation}
Again, we will show that $\bu = \bv$. Suppose not. Then there exist a (possibly empty) 
word $w$ and two infinite words $\bu'$ and $\bv'$ over $\cA$ such that $\bu = wz_1\bu'$ 
and $\bv = wz_2\bv'$ for some letters $z_1$ and $z_2$ with $z_1 \ne z_2$. Hence 
$x\bu = xwz_1\bu'$ and $y\bv = ywz_2\bv'$, and therefore the words $xwz_1$ and $ywz_2$ 
are both factors of $\bt$, since any factor of $x\bu$ and of $y\bv$ is also a factor of 
$\bt$ (by \eqref{eq:inf-1*} and \eqref{eq:inf-2*}). But then, by Richomme's local balance 
property (Theorem~\ref{T:gR07aloc}), $z_2 = x$ or $z_1 = y$. 

If $z_2 = x$, then for any acceptable pair $(x,<_x)$, we have $x <_x z_1$ (since $z_1 \ne z_2$), 
and hence $x\bv ~(= xwx\bv') <_x x\bu ~(= xwz_1\bu')$, contradicting the (lexicographical) 
minimality of $\bu$ with respect to the total order $<_x$. Likewise, if $z_1 = y$, then for 
any acceptable pair $(y,<_y)$, we have $y <_y z_2$ (since $z_1 \ne z_2$), and hence 
$y\bu ~(= ywz_1\bu') <_y  y\bv ~(= ywz_2\bv')$, a contradiction. Thus $\bu = \bv$. 

Hence, there exists a (unique) infinite word $\bu$ on $\cA$ such that, for any acceptable 
pair $(a,<)$, $a\bu \leq T^i(\bt)$ for all $i \geq 0$. 

Conversely, suppose there exists an infinite word $\bu$ on $\cA$ such that, for any 
acceptable pair $(a,<)$, we have
\begin{equation} \label{eq:epi-inequalities}
  a\bu \leq T^i(\bt) \quad \mbox{for all $i \geq 0$}.
\end{equation}
Let $z$ be a left special factor (if any) of $\bt$, and let $z'$ denote the prefix of 
$\bu$ with $|z'| = |z|$. Since $z$ is left special in $\bt$, there exist at least two 
distinct letters $x$, $y$ such that $xz$ and $yz$ are both factors of $\bt$. In particular, 
there exist non-negative integers $\ell_1$ and $\ell_2$ such that $T^{\ell_1}(\bt)$ begins 
with $xz$ and $T^{\ell_2}(\bt)$ begins with $yz$. 
Thus, by inequality~\eqref{eq:epi-inequalities}, we have
\[
xz' \leq_x xz \quad \mbox{for any acceptable pair $(x,<_x)$},
\]
and
\[
yz' \leq_y yz \quad \mbox{for any acceptable pair $(y,<_y)$}.
\]
Hence $z' \leq_x z$ and $z' \leq_y z$, and this implies that $z = z'$. Therefore $\bt$ has 
at most one left special factor of each length and the left special factors of $\bt$ are 
exactly the prefixes of $\bu$. Thus $F(\bu) \subseteq F(\bt)$; in particular, $\bu$ is in 
the closure of the shift-orbit of $\bt$. 

Now suppose that $\bt$ is not episturmian. Then, by Theorem~\ref{T:gR07aloc}, there exists 
a word $w$ (possibly empty) and letters $a$, $b$, $c$, and $d$ with $\{a,b\} \cap \{c,d\} 
 = \emptyset$ such that $awb$ and $cwd$ are both factors of $\bt$. Since $a \ne c$, the word 
$w$ is a left special factor of $\bt$, and therefore $w$ is a prefix of $\bu$. 

Let $\ell_1$ and $\ell_2$ be non-negative integers such that $T^{\ell_1}(\bt)$ begins 
with $awb$ and $T^{\ell_2}(\bt)$ begins with $cwd$. Then, for any two acceptable pairs 
$(a,<_a)$ and $(c,<_c)$, we have
\begin{equation} \label{eq:epi-1}
a\bu~(= awz \cdots) \leq_a T^{\ell_1}(\bt) ~(= awb\cdots),
\end{equation}
and
\begin{equation} \label{eq:epi-2}
c\bu~(= cwz \cdots) \leq_c T^{\ell_2}(\bt) ~(= cwd\cdots).
\end{equation}
Inequality~\eqref{eq:epi-1} implies that $z \leq_a b$, whereas inequality~\eqref{eq:epi-2} 
implies that $z \leq_c d$, and moreover $z \leq_c b$ and $z \leq_a d$. These inequalities 
imply that $z = b = d$, a contradiction.

Hence $\bt$ is episturmian, and therefore $\bu$ is episturmian too (since $\bu$ is in the 
closure of the shift-orbit of $\bt$, which consists of all episturmian words with the same 
set of factors as $\bt$ -- see Remark~\ref{R:episturmian-orbit} or \cite{aGjJ07epis}). 
Moreover, $\bu$ is epistandard since all of its left special factors are prefixes of it. 
Therefore, for any letter $x$ in $\cA$, $x\bu$ is episturmian if and only if $x$ occurs 
infinitely often in the directive word of $\bu$ (see \cite[Theorem 3.17]{jJgP02epis}, 
\cite[Theorem~2.6]{aG07orde}, or \cite[Theorem~6]{gR07aloc}). Hence, for any acceptable 
pair $(a,<)$, $a\bu = \inf\{T^k(\bt), k\geq 0\}$ if and only if the letter $a$ occurs 
infinitely often in the directive word of $\bu$. 
\qed
\end{proof}


\begin{remark}
An unrelated connection between finite balanced words (i.e., finite Sturmian words) 
and lexicographic ordering was recently studied by Jenkinson and Zamboni \cite{oJlZ04char},
who presented three new characterizations of ``cyclically'' balanced finite words via 
orderings. Their characterizations are based on the ordering of shift-orbits, 
either lexicographically or with respect to the $1$-{\em norm} $\mid\cdot\mid_1$, which 
counts the number of occurrences of the symbol $1$ in a given finite word over $\{0,1\}$. 
\end{remark}

\subsection{Sturmian morphisms}\label{stumorp}

Prior to the recent work of Pirillo and others, the extremal property
\eqref{eq:characteristic} was shown to hold for characteristic Sturmian   
sequences in a paper by Berstel and S{\ee}bold \cite{jBpS93acha}. Here is a
reformulation of their result
(recalling the definition of $s_{\alpha,\rho}$ from Section \ref{SS:Sturmian},
and letting $\bc_\alpha := s_{\alpha,\alpha} = s'_{\alpha,\alpha}$ denote
the unique characteristic Sturmian sequence of slope $\alpha$):

\begin{proposition} \label{P:jBpS93acha-1} {\rm \cite[Property 7]{jBpS93acha}}
Let $\alpha > 0$ be an irrational number. Then, for all $i\geq 1$, we have
\[
 a\bc_\alpha < T^i(a\bc_\alpha) \quad \mbox{and} \quad b\bc_\alpha > T^i(b\bc_\alpha).
\]
In particular, for all $i\geq 0$, we have
\[
 a\bc_\alpha < T^i(\bc_\alpha) < b\bc_\alpha.
\]
\end{proposition}

\begin{remark} \label{R:orbit-closure}
Recall from Remark~\ref{R:episturmian-orbit} that the closure of the
shift-orbit of any Sturmian word $\bs$ is a minimal dynamical system
consisting of all the Sturmian words with the same set of factors as
$\bs$ (also see \cite[Proposition~2.1.25]{mL02alge}). In particular,
if $\bs$ is a Sturmian word with (irrational) slope $\alpha$, then
$\Closure(\bs)$ consists of all Sturmian words of slope $\alpha$
(e.g., see \cite[Propositions 2.1.18]{mL02alge} or \cite{fM89infi}).
Accordingly, the second part of Proposition~\ref{P:jBpS93acha-1} (also
see Theorems~\ref{P:JPA} and \ref{T:intro}) tells us that $a\bc_\alpha$ and
$b\bc_\alpha$ are the lexicographically least and greatest Sturmian
words of slope $\alpha$, respectively.
\end{remark}

Proposition \ref{P:jBpS93acha-1} was also proved by Borel and Laubie \cite{jBfL93quel}
in the same year (1993). In \cite{jBpS93acha}, Berstel and S{\ee}bold showed that it is
an easy consequence of the following more general result.

\begin{proposition} \label{P:jBpS93acha-2} 
Let $\alpha > 0$ be an irrational number and let $\rho$, $\rho'$ be real numbers
such that $0 \leq \rho$, $\rho' < 1$. Then
\[
  s_{\alpha,\rho} < s_{\alpha,\rho'} \quad \Longleftrightarrow \quad \rho < \rho'.
\]
\end{proposition}

The above proposition was one of numerous results in \cite{jBpS93acha} leading to the proof 
of a now well-known characterization of {\em Sturmian morphisms}, i.e., morphisms that preserve 
Sturmian words. Specifically, a morphism on $\{a,b\}$ is Sturmian if and only if it can be 
expressed as a finite composition of the following morphisms, in any number and order:
\[
  E: \begin{matrix}
      &a &\mapsto &b& \\
      &b &\mapsto &a&
     \end{matrix}, \qquad \varphi: \begin{matrix}
                                   &a &\mapsto &ab& \\
                                   &b &\mapsto &a&
                                   \end{matrix}, \qquad
  \rev{\varphi}: \begin{matrix}
                  &a &\mapsto &ba& \\
                  &b &\mapsto &a~&
                 \end{matrix}.
\]
(Note that $\varphi = \psi_a\theta_{ab}$ and
$\tilde\varphi = \bar\psi_a\theta_{ab}$; see Section~\ref{SS:episturmian}.)

This result played a particularly important role in Berstel and S{\ee}bold's
characterization of morphisms that preserve characteristic Sturmian words --
the so-called {\em characteristic} or {\em standard (Sturmian) morphisms}.
That is, a morphism on $\{a,b\}$ is {\em standard} if and only if it is
expressible as a finite composition of the morphisms $E$ and $\varphi$ in any
number and order \cite{jBpS93acha}. The fact that there is no occurrence of
the morphism $\rev{\varphi}$ in such a composition is due to Proposition~\ref{P:jBpS93acha-1}.

\subsection{The lexicographic world}

As mentioned previously, a disguised form of Theorem~\ref{T:intro} 
(also see \eqref{eq:characteristic}) appeared in S.~Gan's paper~\cite{sG01stur}; 
in fact, as we shall see, Theorem~\ref{P:JPA} can be deduced from the main results 
in \cite{sG01stur}. Gan came across this  property of Sturmian sequences whilst 
endeavouring to obtain a complete description of the {\em lexicographic world}, defined as follows.

For any two infinite words $\bx$, $\by \in \{0,1\}^\omega$, define the set
\[
\varSigma_{\bx\by} := \{\bs \in \{0,1\}^\omega, \,  \forall i \geq 0, \, \bx \leq T^i(\bs) \leq 
\by\}. 
\] 
The {\em lexicographic world} $\mathcal{L}$ is defined by
$$
\mathcal{L} := \{(\mathbf{x},\mathbf{y}) \in \{0,1\}^\omega \times
\{0,1\}^\omega, \, \varSigma_{\mathbf{xy}} \neq \emptyset\}.
$$
Gan proved in \cite[Lemma~2.1]{sG01stur} that
$$
\cL = \{(\bu,\bv) \in \{0,1\}^\omega \times \{0,1\}^\omega, \, \bv \geq \phi(\bu)\},
$$ 
where $\phi : \{0,1\}^\omega \rightarrow \{0,1\}^\omega$ be the map defined by 
\[
  \phi(\bx) := \inf\{\by \in \{0,1\}^\omega, \, \varSigma_{\bx\by} \ne \emptyset\}.
\]
As Gan points out in that paper, the set $\cL$ is closely related to the bifurcation 
of a Lorenz-like map (see \cite{rLsP01bifu} for example). 

The following theorem combines Corollary~5.6 and Theorem~5.7 from Gan's paper~\cite{sG01stur} 
(also see Theorem~1.1 in the same paper).  It shows in particular that any element in the image 
of $\phi$ is a Sturmian or periodic balanced sequence in $\{0,1\}^\omega$ (and such sequences are the lexicographically greatest amongst their shifts). 

\begin{theorem} \label{T:Gan1} For any sequence $\bs \in \{0,1\}^\omega$, the following 
conditions are equivalent.
\begin{itemize}
\item[(i)] $\bs = \phi(\bx)$ for some sequence $\bx \in \{0,1\}^\omega.$
\item[(ii)] $\bs$ is a Sturmian or periodic balanced sequence satisfying $T^i(\bs) \leq \bs$ 
for all $i \geq 0$.
\end{itemize}
Moreover, if $\bx$ begins with $1$, then $\phi(\bx) = 1^\omega$, and if $\bx = 0\bu$ for some 
$\bu \in \{0,1\}^\omega$, then $\phi(\bx)$ is the unique Sturmian or periodic balanced sequence 
$\bs$ in $\{0,1\}^\omega$ satisfying $0\bu \leq T^i(\bs) \leq 1\bu$ and $T^i(\bs) \leq \bs$ for all 
$i \geq 0$. 
\end{theorem}

In the process of establishing Theorem~\ref{T:Gan1}, Gan also proved the following description 
of {\em Sturmian minimal sets} (see \cite{gH44stur} for a definition; also note that minimal sets 
correspond to minimal dynamical systems).

\begin{theorem} \label{T:Gan2} {\em \cite{sG01stur}}
A minimal set $M$ is a Sturmian minimal set if and only if 
$M \subseteq [0\bx, 1\bx] := \{\by \in \{0,1\}^\omega, 0\bx \leq \by \leq 1\bx\}$ for some 
$\bx \in \{0,1\}^\omega$. Moreover, for any $\bx \in \{0,1\}^\omega$, there exists a unique 
Sturmian minimal set in $[0\bx, 1\bx]$.
\end{theorem}

Theorem~\ref{T:Gan2} actually encompasses the first part of Theorem~\ref{P:JPA}; indeed, 
it can be interpreted as follows: a uniformly recurrent sequence $\by \in \{0,1\}^\omega$ 
satisfies $0\bx \leq T^i(\by) \leq 1\bx$ for all $i\geq 0$ and some binary sequence $\bx$  
if and only if $\by$ is a Sturmian or periodic balanced sequence. As discussed in 
Section~\ref{SS:Pirillo}, this result was recently rediscovered by Glen, Justin, and 
Pirillo~\cite{aGjJgP06char} (see~\eqref{e1}), but in a slightly stronger form without the 
uniform recurrence condition, giving that $\by$ is either a Sturmian sequence, a periodic 
balanced sequence, or a skew sequence (i.e., $\by$ is a balanced sequence).

The second part of Theorem~\ref{P:JPA} can also be deduced from Gan's work, as follows. 
Let $\bu$ be any characteristic Sturmian sequence on $\{0,1\}$. Then, by Theorem~\ref{T:Gan1}, 
the sequence $\bs := \phi(0\bu)$ is the unique Sturmian sequence satisfying 
$0\bu \leq T^i(\bs) \leq 1\bu$ and $T^i(\bs) \leq \bs$ for all $i \geq 0$. 
Suppose $\bx$ is the unique characteristic Sturmian sequence in $\Closure(\bs)$, 
the closure of the shift-orbit of $\bs$. Then $0\bx$ and $1\bx$ are Sturmian sequences, 
by \cite[Proposition 2.1.22]{mL02alge}. Moreover, $0\bx$ and $1\bx$ have the 
same set of factors as $\bx$ since the prefixes of $\bx$ are exactly its left special 
factors. Hence, both $0\bx$ and $1\bx$ are in $\Closure(\bs)$, and therefore, since 
$0\bu \leq T^i(\bs) \leq 1\bu$ for all $i \geq 0$, we have $0\bu \leq 0\bx$ and $1\bx \leq 1\bu$. 
These inequalities imply that $\bu = \bx$. Thus, for any characteristic Sturmian sequence $\bx$, 
we have $0\bx < T^i(\bx) < 1\bx$ for all $i \geq 0$. This establishes the forward direction of  
Theorem~\ref{T:intro}, and it follows that for any Sturmian sequence $\bs$ on $\{0,1\}$, we 
have $0\bu \leq T^i(\bs) \leq 1\bu$ for all $i \geq 0$, where $\bu$ is the unique characteristic 
Sturmian sequence with the same slope as $\bs$ (recall Remark~\ref{R:orbit-closure}). This proves 
the second part of Theorem~\ref{P:JPA} and from this theorem one can easily deduce both 
directions of Theorem~\ref{T:intro} (see Remark~\ref{2.1-2.2}).

\begin{remark} By Remark~\ref{R:orbit-closure}, the lexicographically greatest and least 
Sturmian sequences in the closure of the shift-orbit of a Sturmian sequence $\bs$ on $\{0,1\}$ 
are $0\bu$ and $1\bu$ where $\bu$ is the unique characteristic Sturmian sequence with the same 
slope as $\bs$. We thus deduce from Theorems~\ref{P:JPA} and \ref{T:Gan1} that, for any 
sequence $\bx$ on $\{0,1\}$ beginning with $0$, the sequence $\phi(\bx)$ is a Sturmian or 
periodic balanced sequence of the form $1\bu$. Moreover, if $\phi(\bx)$ is Sturmian, then 
$\bu$ is the unique characteristic Sturmian sequence with the same slope as $\phi(\bx)$.
\end{remark}

The following lemma was a key step in Gan's proof of Theorem \ref{T:Gan2}. It involves the 
{\em block condition} (BC): a sequence $\bs \in \{0,1\}^\omega$ satisfies the BC if, for any 
finite word $w$ on $\{0,1\}$, at least one of the words $0w0$ and $1w1$ is not a factor of $\bs$. 

\begin{lemma} \label{L:Gan4.4} {\rm \cite[Lemma 4.4]{sG01stur}}
A sequence $\bs \in \{0,1\}^\omega$ satisfies the BC if and only if there exists a 
sequence $\bu$ such that $0\bu \leq T^i(\bs) \leq 1\bu$ for all $i\geq 0$.
\end{lemma}

This result is essentially the characterization of balanced infinite words given in 
\cite{aGjJgP06char} (see \eqref{e1}). Indeed, the BC is equivalent to the balance 
property, as defined in Definition \ref{D:balance}. See Section~3 in \cite{eCgH73sequ}, 
in which the balance property is called the {\em Sturmian block condition} (see also 
\cite{gR07aloc}). Note that the BC of Coven and Hedlund 
\cite[Lemma 3.06 p.\ 143]{eCgH73sequ} is stronger than Gan's in that 
``for any finite word $w$'' is replaced by ``for any palindrome $w$''.

\begin{remark} As explained by Labarca and Moreira in \cite{LabMor0}, the terminology 
``lexicographical world'' was coined in 2000, in a preprint version of \cite{LabMor2} 
(which appeared only in 2006) in which the authors extended the work of Hubbard and
Sparrow \cite{jHcS90thec}. For more on the lexicographic(al) world, the reader can 
look at, e.g., \cite{LabMor1, LabMor2} and the references therein. See also the recent paper \cite{jAaG09dist}, in which the present two authors give a complete description of the lexicographic world in the process of describing the minimal intervals containing all fractional parts $\{\xi2^n\}$, for some positive real number $\xi$, and for all $n\geq 0$.
\end{remark}

\subsection{The early work of Veerman: 1986 \& 1987}\label{veer}

Let $\cS^\alpha$ denote the set of all Sturmian sequences of (irrational)
slope $\alpha > 0$ over the alphabet $\{0,1\}$ (i.e., $a\mapsto
0$, $b \mapsto 1$ in Theorem~\ref{stu}). As noted, e.g., in
\cite{jB02rece}, each Sturmian sequence $\bs \in \cS^\alpha$ can be viewed
as the binary expansion of some real number $r(\bs)$ modulo $1$. Moreover,
it is easily verified that, for any $\bs$, $\bs' \in \cS^\alpha$, we have
$\bs < \bs'$  if and only if $r(\bs) < r(\bs')$.  Furthermore, by
Remark~\ref{R:orbit-closure}, we know that the lexicographically least and
greatest sequences in $\cS^\alpha$ are $0\bc_\alpha$ and $1\bc_\alpha$,
respectively. In terms of binary expansions, as $r(1\bc_\alpha) = 
1/2 + r(0\bc_\alpha)$, it follows that the set $r(\cS^\alpha) := 
\{r(\bs) \in [0,1),~\bs \in \cS^\alpha\}$ is completely contained within the 
closed interval $[r(0\bc_\alpha), r(1\bc_\alpha)]$ of length $1/2$ and not in any 
smaller interval. This latter result (to compare with Bugeaud-Dubickas' result 
where base $2$ is replaced with base $b$ \cite{yBaD05frac}) is essentially 
a reformulation of Theorem~2 p.~558 in Veerman's paper~\cite{pV86symb}, which 
also states that $r(\cS^\alpha)$ is a Cantor set and that 
$[r(0\bc_\alpha), r(1\bc_\alpha)]$ has
Lebesgue measure zero. The converse of this theorem was proved one year
later by Veerman in \cite[Theorem 2.1, p. 193--194]{pV87symb}. As such, it
seems that Veerman was the first to (implicitly) prove the Sturmian extremal
property given in Theorem~\ref{P:JPA}, under the framework of symbolic
dynamics. Actually, Veerman's main result in \cite{pV87symb} shows that a
sequence $\bs$ in $\{0,1\}^\omega$ satisfies the inequalities $0\bu \leq T^i(\bs)
\leq 1\bu$ for some sequence $\bu \in \{0,1\}^\omega$ and for all $i \geq 0$ 
if and only if $\bs$ is a Sturmian sequence or a periodic balanced sequence 
({\it cf.}~\eqref{e1}). A few years earlier (in 1984),  Gambaudo 
{\it et al.}~\cite{jGoLcT84dyna} had already proved the
periodic case (i.e., the case when $\alpha$ is rational); Veerman considered
his Theorem~2.1 in \cite{pV87symb} to be a generalization of their main
result.

\begin{remark}
Note that the set $r(\cS^\alpha)$ is a dynamical system under the operation
of the {\em doubling map} $\sigma: x \mapsto 2x \pmod{1}$ on the
one-dimensional torus $\TT = \RR/\ZZ$. This was the point of view of 
Veerman and also that of Bullet and Sentenac~\cite{sBpS94orde}, who gave
reformulations and self-contained combinatorial proofs of some of Veerman's
results in~\cite{pV86symb, pV87symb}. In particular, Bullett and Sentenac
gave another proof of the following result (which can be deduced from
Veerman's work): for each closed interval $C_{\mu} = [\mu, 1/2+\mu]$ of
length $1/2$ (where $\mu \in \TT$), there exists a unique $\alpha$ such that
$r(\cS^\alpha)$ is contained in $C_\mu$ and there is no other dynamical
system for the doubling map that is a strict subset of $C_\mu$. This fact was
recently used by Jenkinson~\cite{oJ07opti} to prove new characterizations of
{\em Sturmian measures}, which have applications to ergodic optimization of
convex functions. Another important application is in the combinatorial
description of the Mandelbrot set (e.g., see \cite{sBpS94orde, kK00inva}).
\end{remark}

\begin{remark} In the study of kneading sequences of {\em Lorenz maps}
(i.e., a certain class of piece-wise monotonic maps on $[0,1]$ with a
single discontinuity), Glendinning, Hubbard, and
Sparrow~\cite{pGcS93prim, jHcS90thec} have investigated so-called 
{\em allowed pairs} $(\br,\bs)$ of distinct binary sequences in 
$\{0,1\}^\omega$ satisfying
\[
\br \leq T^i(\br) < \bs \quad \mbox{and} \quad 
\br < T^i(\bs) \leq \bs \quad \mbox{for all $i\geq 0$}.
\]
In particular, it was shown in \cite{jHcS90thec} that these allowed 
pairs are exactly the pairs of (distinct) binary sequences in 
$\{0,1\}^\omega$ that are realizable as kneading invariants 
of a topologically expansive Lorenz map. (Note that the case 
$\bs = 1^\omega$ was studied by Acquier, Cosnard, and Masse in 
\cite{ACM}.) Moreover it can be deduced from
property~\eqref{eq:characteristic} that the allowed pairs of the form 
$(0\bu, 1\bu)$ are those where $\bu$ is a characteristic Sturmian sequence.
\end{remark}

\section{Back to distribution modulo $1$: 
the Thue-Morse sequence shows up}\label{morsesection}

As indicated in the introduction, we began writing this survey after the publication 
of the paper of Bugeaud and Dubickas \cite{yBaD05frac}, whose starting point goes
back to a paper of Mahler \cite{kM68anun}. In that paper Mahler defines the set of
$Z$-numbers 
$$
\left\{
\xi \in {\mathbb R}, \ \xi > 0, \ \forall n \geq 0, \ 
0 \leq \left\{ \xi \left(\frac{3}{2}\right)^n\right\} < \frac{1}{2} 
\right\}
$$
where $\{x\}$ is the fractional part of the real number $x$.
Mahler proved that this set is at most countable. It is still an open problem to 
prove that this set is actually empty. More generally, given a real number $\alpha > 1$
and an interval $(s,t) \subset (0,1)$ one can ask whether there exists $\xi > 0$
such that, for all $n \geq 0$, we have $s \leq \{\xi \alpha^n\} < t$.
Flatto, Lagarias, and Pollington \cite[Theorem~1.4]{FLP} proved that, if $\alpha = p/q$ 
with $p, q$ coprime integers and $p > q \geq 2$, then any interval $(s, t)$ such that
for some $\xi > 0$, one has that $\{\xi (p/q)^n\} \in (s, t)$ for all $n \geq 0$,
must satisfy $t-s \geq 1/p$. The main result in \cite{yBaD05frac} reads as follows.

\begin{theorem}[Bugeaud-Dubickas]\label{bugdub}
Let $b \geq 2$ be an integer and let $\xi$ be an irrational number. Then the numbers
$\{\xi b^n\}$ cannot all lie in an interval of length $< 1/b$. Furthermore there
exists a closed interval $I$ of length $1/b$ containing the numbers $\{\xi b^n\}$
for all $n \geq 0$ if and only if the sequence of base $b$-digits of the fractional 
part of $\xi$ is a Sturmian sequence $\bs$ on the alphabet $\{k, k+1\}$ for some 
$k \in \{0, 1, \ldots, b-2\}$. If this is the case, then $\xi$ is transcendental,
and the interval $I$ is semi-open. It is open unless there exists an integer $j \geq 1$ 
such that $T^j(\bs)$ is a characteristic Sturmian sequence on the alphabet 
$\{k, k+1\}$.
\end{theorem}

The reader will easily see the relation between Theorem~\ref{bugdub} and
Theorems~\ref{P:JPA} and \ref{T:intro}. 
Note that the first assertion in Theorem~\ref{bugdub} is generalized to
algebraic real numbers $> 1$ by Dubickas in \cite{aD06arith}. Also note that 
two other papers by Dubickas \cite{aD06onth, aD07onas} deal with links between 
distribution of $\{\xi \alpha^n\}$ modulo $1$ and combinatorics on words. 
Furthermore the {\em Thue-Morse sequence}, defined as the fixed point beginning
with $0$ of the morphism $0 \to 01$, $1 \to 10$, shows up in these two papers:
in \cite{aD06onth} for the study of ``small'' and ``large'' limit 
points of $\Vert \xi(p/q)^n\Vert$, the distance to the nearest 
integer of the product of any non-zero real number $\xi$ by the 
powers of a rational; in \cite{aD07onas} for the study of  the 
``small'' and ``large'' limit points of the sequence of fractional 
parts $\{\xi b^n\}$, where $b < -1$ is a negative rational number 
and $\xi$ is a real number. For work in a similar vein and with an avatar of
the Thue-Morse sequence, see \cite{Kaneko}.

\medskip

Interestingly enough, the Thue-Morse sequence also appeared in 1983 in 
another question of distribution, as a by-product of the combinatorial 
study of a set of sequences related to iterating continuous maps 
of the unit interval (see \cite{Allouche, AllCos1}). 

\begin{theorem}\label{etat}
Define the set $\widetilde{\Gamma}$ by
$$
\widetilde{\Gamma} := \{x \in [0, 1], \ 1-x \leq \{2^kx\} \leq x \}.
$$
Then the smallest limit point of $\widetilde{\Gamma}$ is the number 
$\alpha:=\sum a_n/2^n$, where $(a_n)_{n \geq 0}$ is the Thue-Morse 
sequence. The set $\widetilde{\Gamma}$ contains only countably many 
elements less than $\alpha$ and they are all rational. Furthermore any 
segment on the right of $\alpha$ contains uncountably many elements of 
$\widetilde{\Gamma}$. This structure around $\alpha$ repeats at 
infinitely many scales: $\widetilde{\Gamma}$ is a fractal set.
\end{theorem}

The reader will have guessed that Theorem~\ref{etat} above is a by-product
of the combinatorial study of the set
\begin{equation} \label{eq:Gamma}
\Gamma := \{\bu \in \{0, 1\}^{\mathbb N}, \ \forall k \geq 0, \ 
\overline{\bu} \leq T^k(\bu) \leq \bu\}
\end{equation}
where $\overline{\bu}$ is the sequence obtained by switching $0$'s 
and $1$'s in $\bu$ (see \cite{Allouche}).

\bigskip

An avatar of the set $\Gamma$ (where large inequalities are replaced by strict
inequalities) was studied in \cite{ErdJooKom} in the description of {\em univoque
numbers}, i.e., real numbers $\beta$ in $(1, 2)$ such that there exists a unique
base $\beta$-expansion of $1$ as $1 = \sum_{j \geq 1} u_j \beta^{-j}$, with 
$u_j \in \{0, 1\}$. See \cite{AllCos2} for more details.

In \cite{jA07anot} JPA uses Theorem~\ref{P:JPA} to prove that a Sturmian sequence 
$\bs$ on $\{0,1\}$ belongs to the set $\Gamma$ (see \eqref{eq:Gamma}) if and only if there 
exists a characteristic Sturmian sequence $\bu$ beginning with $1$ such that $\bs = 1\bu$. 
(In particular, a Sturmian sequence belonging to $\Gamma$ must begin with $11$.) 
As an immediate corollary we have that a real number $\beta \in (1,2)$ is univoque
and {\em self-Sturmian} (i.e., the greedy $\beta$-expansion of $1$ is 
a Sturmian sequence) if and only if the $\beta$-expansion of $1$ is of the form $1\bu$, where 
$\bu$ is a characteristic Sturmian sequence beginning with $1$.
Self-Sturmian numbers were introduced in \cite{dCdK04stur}, where it was proved that such
numbers are transcendental (see also \cite{kwon} for more on related questions). 
Theorem~\ref{T:intro} was used in \cite{dCdK04stur} and a proof of Theorem~\ref{P:JPA} was 
also given in a preprint version of that paper (see {\tt http://arxiv.org/abs/math/0308140}); 
it was taken off the last version, as the author explained to JPA: first because a referee 
suggested it was ``folklore'', and second because actually only one direction of 
Theorem~\ref{T:intro} was needed. Self-sturmian numbers have since been generalized to 
{\em self-episturmian numbers} in \cite{aG07orde}, where an analogue of Theorem~\ref{P:JPA} 
for episturmian sequences can also be found (see Theorem~\ref{AG}).

Also note that sets related to the set $\Gamma$ and to the lexicographic world
occur in the study of badly approximable numbers in \cite{Nilsson}.

\bigskip

We end this section with a last remark.

\begin{remark}
It is tempting to try to convert the extremal property for episturmian sequences
given in Corollary~\ref{Cor:aGjJgP06char-infinite} (see \cite{aGjJgP06char}) to
a result in distribution modulo $1$. From now on, $<$ will denote the ``usual''
order on $D:= \{0, 1, \ldots, d-1\}$; other orders will be denoted by $\prec$.
As we have seen, an infinite word $\bt$ on $D:= \{0, 1, \ldots, d-1\}$ is episturmian 
in the wide sense (i.e., episturmian or episkew) if and only if there exists an infinite 
word $\bu$ such that
\begin{equation*}
\hskip\hsize minus 1fill 
a\bu \preceq \min(\bt) 
\hskip\hsize minus 1fill (*)
\end{equation*}
for any acceptable pair $(a, \prec)$.
Actually, replacing the ``usual'' order on $D$ by another total order is 
the same as keeping the order but replacing each $j$ in this set by $\sigma(j)$, 
where $\sigma$ is a permutation of $D$.
More precisely, $(a, \prec)$ is an acceptable pair if and only if there
exists a permutation $\sigma_{\preceq}$ of $D$ such that $\sigma(a) = 0$
and $i \preceq j \Leftrightarrow \sigma(i) \leq \sigma(j)$.
Hence, another way of formulating $(*)$ above is as follows:
there exists an infinite word $\bu$ such that for all permutations $\sigma$
of $D$ one has
\begin{equation*}
0\sigma(\bu) \le \min(\sigma(\bt))
\end{equation*}
where $\sigma(u_0 u_1 u_2 \ldots) := \sigma(u_0) \sigma(u_1) \sigma(u_2) \ldots$
(for finite or infinite words on $D$).
Hence translating extremal properties of episturmian sequences to properties of
distribution modulo $1$ for real numbers consists of looking at reals $x$ in $(0,1)$
such that there exists a real $y$ in $(0,1)$ with 
$\frac{1}{2}y_{\sigma} \leq \{d^k x_{\sigma}\}$ for all integers $k$ and for all 
permutations $\sigma$ (where $x_{\sigma}$ is the real number obtained from $\sigma$ 
by applying the permutation $\sigma$ digitwise). If $d=2$, permuting $0$'s and $1$'s 
in a real number $x$ written in base $2$ is the same as replacing $x$ by $1-x$. 
Hence, in that case, the inequalities $\frac{1}{2}y_{\sigma} \leq \{2^k x_{\sigma}\}$ 
boil down to the two families of inequalities $\frac{1}{2}y \leq \{2^k x\}$ and 
$\frac{1}{2}(1-y) \leq \{2^k(1-x)\} = 1 - \{2^k x\}$, i.e., 
$\frac{1}{2}y \leq \{2^k x\} \leq \frac{1}{2} + \frac{1}{2}y$ for all $k$.
This is precisely the question from which we started our paper, but for general $d$ 
it does not seem that number-theoretists have been interested in distribution modulo $1$
combined with permuting digits.
\end{remark}

\section{Addendum}

While writing this survey we came across several extra references that are 
related to some of its parts. We give some of them here; the interested reader can
look at these papers and the references therein: about combinatorics of words 
and Lorenz maps \cite{AlFa1, AlFa2, AlMa1, AlMa2, Keller-StP, SiSo}, about extremal 
properties of Sturmian sequences or measures \cite{Jenkinson08, Kieffer, KSWZ}, about 
the distribution of $\{\xi\alpha^n\}$ \cite{Aki, AFS, Dub2009, Zaimi1, Zaimi2}, and 
last but not least the historical paper of Lorenz \cite{Lorenz} (also see \cite{Sparrow}).

\section{Acknowledgements} 

The authors would like to thank J. Berstel, J. Cassaigne, J. Justin, D. Kwon,
G. Pirillo, G. Richomme, P. S\'e\'ebold, and L. Q. Zamboni for discussions, comments, 
and suggestions.

\end{document}